\documentclass{article}

\usepackage[T2A]{fontenc}
\usepackage[cp1251]{inputenc}
\usepackage[english,russian]{babel}
\usepackage[tbtags]{amsmath}
\usepackage{amsfonts,amssymb}

\usepackage{mathrsfs}

\usepackage{graphicx}

\numberwithin{equation}{subsection}

\newtheorem{theorem}{Теорема}

\newcommand{\Natur}{{\mathbb N}}

\newtheorem{claim}{Утверждение}

\begin{document}


\date{24.06.2015}

\author{Ф.\,А.~Пушняков}

\title{О числе ребер в индуцированных подграфах специального дистанционного графа}

\markboth{Ф.\,А.~Пушняков}{О числе ребер в индуцированных подграфах специального дистанционного графа}

\maketitle

\begin{abstract}
В работе получены новые оценки числа ребер в индуцированных подграфах специального дистанционного графа.
Библиография: 21 название.
\end{abstract}


\section{Введение}

Рассмотрим последовательность графов $G_{n} = G_{n} (V_{n}, E_{n}) = G(n,3,1)$,  у которых $$V_{n} = \{x = (x_{1}, \dots, x_{n}) \; \vert \; x_{i} \in \{0, 1\},\; i = 1, \dots , n \; , \; x_{1} + \ldots + x_{n} = 3\},$$ $$E_{n} = \{ (x, y) \; \vert \; \langle x, y \rangle = 1\},$$
где через $\langle x, y \rangle$ обозначено скалярное произведение векторов $x$ и $y$. 
Иными словами, вершинами графа $G(n, 3, 1)$ являются $(0,1)$-векторы, скалярный квадрат которых равен трем. И эти вершины соединены ребром тогда и только тогда, 
когда скалярное произведение соответствующих веторов равно единице. Данное определение можно переформулировать в комбинаторных терминах. 
А именно, рассмотрим граф, вершинами которого являются всевозможные трехэлементные подмножества множества $\mathcal{R}_{n} = \{1, \dots, n\}$, 
причем ребро между такими вершинами проводится тогда и только тогда, когда соответствующие трехэлементные подмножества имеют ровно один общий элемент. 
Изучение данного графа обусловлено многими задачами комбинаторной геометрии, экстремальной комбинаторики, теории кодирования: 
например, задачей Нелсона--Эрдёша--Хадвигера о раскраске метрического пространства (см. \cite{Rai3}--\cite{KW}), 
проблемой Борсука о разбиении множеств в пространствах на части меньшего диаметра (см. \cite{Rai3}--\cite{Rai5}, \cite{Bolt}--\cite{Rai2}), 
задачами о числах Рамсея (см. \cite{Ram}, \cite{Nagy}), задачами о кодах с одним запрещенным расстоянием (см. \cite{MS}, \cite{Bass}).
\par Напомним несколько свойств данного графа. Граф $G(n,3,1)$ является регулярным со степенью вершины $d_{n} = 3 \cdot C_{n-3}^{2}$. Очевидно, что $\displaystyle \vert V_{n} \vert = C_{n}^{3} \sim \frac{n^{3}}{6}$ при  $n \rightarrow \infty$. В силу регулярности рассматриваемого графа имеем $\vert E_{n} \vert = \frac{d_{n} \cdot \vert V_{n} \vert}{2} = \frac{3}{2} \cdot C_{n-3}^{2} \cdot C_{n}^{3} \sim \frac{n^{5}}{8}$ при $n \rightarrow \infty$. \par
Напомним, что \textit{независимым множеством} графа называется такое подмножество его вершин, что никакие две вершины подмножества не соединены ребром. 
\textit{Числом независимости} $\alpha(G)$ называется наибольшая мощность независимого множества.  Положим $\alpha_{n} = \alpha(G(n, 3, 1))$. 
Результат теоремы Ж. Надя (см. \cite{Nagy}) отвечает на вопрос о числе независимости графа $G(n, 3, 1)$. А именно, $\alpha_{n} \sim n$ при $n \rightarrow \infty$. Более того, из доказательства теоремы Ж. Надя можно сделать вывод о структуре независимого множества в рассматриваемом графе. Для описания этой структуры введем дополнительные обозначения.  Пусть $W \subseteq V_{n}$. Будем говорить, что $W$ является множеством вершин \textit{первого типа}, если $\vert W \vert \geq 3$ и существуют такие $i, j \in \mathcal{R}_{n}$, что для любой вершины $w \in W$ выполнено $i, j \in w$; далее, $W$ является множеством вершин \textit{второго типа}, если $\vert W \vert \geq 2$ и существуют такие $i, \; j, \; k, \; t \in \mathcal{R}_{n}$, что для любой вершины $w \in W$ выполнено $w \subset \{i, j, k, t\} $; наконец, $W$ является множеством вершин \textit{третьего типа}, если для любых $w_{1}, w_{2} \in W$ выполнено соотношение $w_{1} \cap w_{2} = \emptyset$. Более того, \textit{носителем} множества вершин назовем объединение всех вершин данного множества. Тогда имеет место следующее утверждение.
\begin{claim} Любое независимое множество $U \subseteq V_{n}$ можно представить в виде объединения $$U = \left( \cup_{i \in \mathcal{I}} A_{i} \right) \cup \left( \cup_{j \in \mathcal{J}} B_{j} \right) \cup \left( \cup_{k \in \mathcal{K}} C_{k} \right),$$ где $A_{i}$ -- множество вершин первого типа, $B_{j}$ -- множество вершин второго типа, $C_{k}$ -- множество вершин третьего типа, $i \in \mathcal{I}, \; j \in \mathcal{J},\;  k \in \mathcal{K}$, и носители всех упомянутых множеств попарно не пересекаются.
\end{claim}
Мы не доказываем данное утверждение, так как оно мгновенно следует из доказательства теоремы Ж. Надя (см. \cite{Nagy}).

\par
Обозначим через $r(W)$ количество ребер графа $G$ на множестве $W \subseteq V_{n}$. Иными словами, $$r(W) = \vert \{(x, y) \in E(G) \; \vert \; x \in W, \; y \in W\} \vert \; .$$ Также положим  $$r(l(n)) = \min_{\vert W \vert = l(n), \; W \subseteq V_{n}} r(W) \; .$$
Заметим, что если $l(n) \leq \alpha_{n}$, то $r(l(n)) = 0$ и обсуждать нечего. Если же $l(n) > \alpha_{n}$, то, очевидно, в любом $W \subseteq V_{n}$ мощности $l(n)$ непременно найдутся ребра. Возникает интересный вопрос об изучении величины $r(l(n))$. В настоящей работе мы приведем практически полное исследование данной величины. Нами доказана следующая теорема.
\begin{theorem}
Имеют место четыре случая:
\begin{enumerate}
\item Пусть функции $f: \Natur \rightarrow \Natur, \; g: \Natur \rightarrow \Natur$ таковы, что выполнено $n = o(f)$ и $g = o(n^{2})$ при $n \rightarrow \infty$. Пусть функция $l: \Natur \rightarrow \Natur$ такова, что для любого $n \in \Natur$  выполнена цепочка неравенств $f(n) \leq l(n) \leq g(n)$. Тогда $r(l(n)) \sim \frac{l(n)^{2}}{2 \alpha_{n}}$ при $n \rightarrow \infty$.
\item Пусть функция $l: \Natur \rightarrow \Natur$  такова, что существуют константы $C_{1}, \; C_{2}$, с которыми для каждого $n \in \Natur$ выполнена цепочка неравенств $C_{1} \cdot n^{2} \leq f(n) \leq C_{2} \cdot n^{2}$. Тогда $r(l(n)) \sim \frac{l(n)^{2}}{2 \alpha_{n}}$ при $n \rightarrow \infty$.
\item Пусть функции $f: \Natur \rightarrow \Natur, g: \Natur \rightarrow \Natur$ таковы, что выполнено $n^{2} = o(f(n))$ и $g(n) = o(n^{3})$ при $n \rightarrow \infty$.  Пусть функция $l: \Natur \rightarrow \Natur$ такова, что для каждого $n \in \Natur$ выполнено $f(n) \leq l(n) \leq g(n)$. Тогда существует такая функция $h: \Natur \rightarrow \Natur$, что $h(n) \sim \frac{5 l(n)^{2}}{\alpha_{n}} $  при $n \rightarrow \infty$ и для каждого $n \in \Natur$ выполнена цепочка неравенств $\frac{l(n)^{2}}{\alpha_{n}} \leq r(l(n)) \leq h(n)$.
\item Пусть функция $l: \Natur \rightarrow \Natur$  такова, что существует константа $C$, с которой выполнена цепочка неравенств $C \cdot n^{3} \leq l(n) \leq  C_{n}^{3}$. Пусть $c_{n} = 1-\frac{l(n)}{C_{n}^{3}}$. Тогда существует функция $f: \Natur \rightarrow \Natur$, такая, что $f(n) \sim n^{5} \left(\frac{1}{8} - \frac{c_{n}}{4} + \frac{c_{n}^{2}}{72}\right)$ при $n \rightarrow \infty$, и для каждого $n \in \Natur$ выполнено $r(l(n)) \geq f(n)$.
\end{enumerate}
\end{theorem}
Проанализируем формулировку данной теоремы. В первых двух случаях мы нашли асимптотическое значение величины $r(l)$ при $n \rightarrow \infty$. В третьем случае мы нашли порядок величины $r(l(n))$. Четвертый случай исследован не до конца, но оценка, полученная в нем, обладает тем свойством, что $r(l(n)) \sim \vert E_{n} \vert $ при $l(n) \sim \vert V_{n} \vert $ и $n \rightarrow \infty$. В следующем разделе мы приведем доказательство теоремы 1.

\section{Доказательство теоремы 1}
\subsection{Доказательство пункта 1}
Нижняя оценка известна и вытекает из классической теоремы Турана (см., например, \cite{Rai11}, \cite{Rai12}). 
Для доказательства верхней оценки необходимо для каждой функции $l(n)$, удовлетворяющей условию пункта 1 теоремы, и для каждого $n$ построить пример множества $W_{n}$ мощности $l(n)$, для которого величина $r(W_{n})$ оценивается сверху нужным образом. При этом можно считать, что $n$ достаточно велико. \par
Зафиксируем произвольную функцию $l$, удовлетворяющую условию пункта 1 теоремы, и число $n$. Положим $a(n) = \left[\frac{n^{2}}{l(n)}\right]$. Положим $b(n) = \left[\ln a(n)\right]$. Положим $x(n) = n - \left[\frac{n}{b(n)}\right]$. Ясно, что $x(n) \sim n$ при $n \rightarrow \infty$. Также положим $y(n) = \left[\frac{2 l(n)}{x(n)}\right]$.
Рассмотрим следующее подмножество множества $\mathcal{R}_{n} = \{1, \dots, n\}$:
$$A_{1} = \{1, \dots, x\}.$$
\par
Рассмотрим также следующее множество вершин:
$$W_{n} = \bigcup_{i \in A_{1}} \bigcup_{j \in \{1, \dots, \left[\frac{y}{2}\right]\}} \{\{x+2(j-1)+1, x+2(j-1)+2, i\}\}.$$
Найдем мощность множества $W_{n}$. Ясно, что $$\vert W_{n} \vert = \vert A_{1} \vert \cdot \left[\frac{y}{2}\right] = x \cdot \left[\frac{y}{2}\right] \sim \frac{xy}{2}$$ при $n \rightarrow \infty$.
Найдем $r(W_{n})$. Обозначим через $E(W_{n})$ множество ребер графа $G(n, 3, 1)$ на множестве вершин $W_{n}$. Иными словами, 
$E(W_{n}) = \{(a, b) \in E(G) \vert \; a \in W_n, \; b \in W_{n}\}$.  \par
Посчитаем мощность множества $E(W_{n})$. Ясно, что только вершины вида $\{x+2(j-1)+1, x+2(j-1)+2, i\}$ при фиксированном $i \in A_{1}$ могут образовывать ребро. Всего существует $\left[\frac{y}{2}\right] \cdot x \cdot (\left[\frac{y}{2}\right] - 1) \cdot \frac{1}{2}$ пар таких вершин. Действительно,  $\left[\frac{y}{2}\right] \cdot x$ способами можно выбрать одну вершину из $W_{n}$,  $(\left[\frac{y}{2}\right] - 1)$ способами можно выбрать ей пару из $W_{n}$, и, наконец, сомножитель $\frac{1}{2}$ показывает нам, что каждую пару вершин мы посчитали два раза. \par
Таким образом, $\vert E(W_{n}) \vert = \left[\frac{y}{2}\right] \cdot x \cdot (\left[\frac{y}{2}\right] - 1) \cdot \frac{1}{2}$. Подставим в полученное выражение значения параметров:
$$\vert E(W_{n}) \vert \sim \frac{x y^{2}}{8} = \frac{x^{2}y^{2}}{8n} \cdot \frac{n}{x} \sim \frac{x^{2}y^{2}}{8n} \sim \frac{l(n)^{2}}{2 \alpha_{n}}.$$
Таким образом, искомая верхняя оценка получена.

\subsection{Доказательство пункта 2}
Нижняя оценка, как и в предыдущем пункте,  известна и вытекает из классической теоремы Турана. Для доказательства верхней оценки необходимо для каждой функции $l(n)$, удовлетворяющей условию пункта 2 теоремы, и для каждого $n$ построить пример множества $W_{n}$ мощности $l(n)$, для которого величина $r(W_{n})$ оценивается сверху нужным образом. По-прежнему можно считать, что $n$ достаточно велико. \par
Зафиксируем произвольную функцию $l$, удовлетворяющую условию пункта 2 теоремы, и число $n$. Положим $c_{n} = 4 - \frac{1}{\ln n}, k = \left[\frac{n}{4}\right]$. Положим
$$W_{1} = \bigcup_{i=3}^{\left[ c_{n}  k \right] } \{\{1,\; 2,\; i\}\},$$
$$W_{2} = \bigcup_{i=3}^{\left[ c_{n}  k \right] } \bigcup_{j=1}^{\left[ \frac{n-\left[ c_{n}k \right]}{2} \right]} \{\{i, \; \left[ c_{n}k \right] + 2(j-1)+1, \; \left[ c_{n}k \right]+2(j-1)+2\}\}.$$
Обозначим $W_{n} = W_{1} \sqcup W_{2}$. Ясно, что $$\vert W_{n} \vert = \vert W_{1} \vert + \vert W_{2} \vert = \left[ c_{n}k \right] - 2 + \left(\left[ c_{n}k \right] - 2\right) \left[ \frac{n-\left[ c_{n}k \right]}{2}\right] \sim c_{n}k\frac{n-c_{n}k}{2} \sim $$ $$ \sim c_{n}k \frac{(4-c_{n})k}{2} = \frac{c_{n}(4-c_{n})k^{2}}{2}.$$ Как и раньше, обозначим через $E(W_{n})$ множество ребер графа $G(n, 3, 1)$ на множестве вершин $W_{n}$. Ясно, что $$\vert E(W_{n}) \vert = (\left[c_{n}k\right]-2)\left[\frac{n-\left[c_{n}k\right]}{2}\right]+\frac{1}{2}(\left[c_{n}k\right]-2)\left[\frac{n-\left[c_{n}k\right]}{2}\right]\left(\left[\frac{n-\left[c_{n}k\right]}{2}\right]-1\right) \sim $$ $$\sim \frac{c_{n} \left(4-c_{n}\right)^{2} k^{3}}{8} = \frac{c_{n}^{2}\left(4-c_{n}\right)^{2}k^{4}}{4}\frac{1}{2c_{n}k} \sim \frac{\vert W_{n} \vert ^{2}}{2 \alpha_{n}}.$$ Таким образом, утверждение пункта 2 доказано.

\subsection{Доказательство пункта 3}
Нижняя оценка вытекает из аналога теоремы Турана для дистанционных графов (см., например, \cite{Rai11}, \cite{Rai12}). Для доказательства верхней оценки необходимо для каждой функции $l(n)$, удовлетворяющей условию пункта 3 теоремы, и для каждого $n$ построить пример множества $W_{n}$ мощности $l(n)$, для которого величина $r(W_{n})$ оценивается сверху нужным образом. По-прежнему можно считать, что $n$ достаточно велико. \par
Зафиксируем произвольную функцию $l$, удовлетворяющую условию пункта 3 теоремы, и число $n$. Положим $k(n) = \left[\frac{l(n)}{\left[ \frac{n}{2}\right] \cdot \left[ \frac{n}{4}\right]}\right]$. Ясно, что $k(n) = o(n)$, при $n \rightarrow \infty$. Рассмотрим следующие подмножества множества $\mathcal{R}_{n}$:
$$A_{1} = \begin{cases}
     \{1, \dots, 2m\}     & \text{при $n=4m$,}\\
  \{1, \dots, 2m+1\}  & \text{при $n=4m+1$,} \\
  \{1, \dots, 2m+2\}  & \text{при $n=4m+2$,} \\
  \{1, \dots, 2m+3\}  & \text{при $n=4m+3$,} \\
  
\end{cases}$$
$$A_{2} = \mathcal{R}_{n} \setminus A_{1}.$$
Ясно, что $\vert A_{1}\vert \sim \frac{n}{2}$, $\vert A_{2}\vert \sim \frac{n}{2}$ при $n \rightarrow \infty$. Также ясно, что число $\vert A_{2} \vert$ четно.
Положим $a(n) = \vert A_{1} \vert$. Пусть $\sigma \in S_{n - a(n)}$ -- произвольная перестановка. Назовем разбиением множества $A_{2}$, отвечающем перестановке $\sigma$, следующее множество:
$$P_{\sigma} = \{\left( a(n) + \sigma(1), \; a(n) + \sigma(2)\right), \dots, \left(a(n) + \sigma(n-a(n)-1), \; a(n)+\sigma(n-a(n))\right)\}.$$
 Ясно, что можно выбрать $k(n)+1$ различных перестановок так, что никакая пара элементов $(x, y) \in A_{2} \times A_{2}$ не будет принадлежать более чем одному разбиению. Иными словами, можно выбрать $k(n)+1$ попарно не пересекающихся разбиений. Обозначим их $P_{1}, \dots, P_{k(n)+1}.$ Тогда для $ i=1, \dots, k(n)+1$ положим
$$W^{(i)} = \bigcup_{x \in A_{1}, \; (y,z) \in P_{i}}  \{\{x, y, z\}\} . $$
Пусть $w(n) = \vert W^{(1)} \vert = \ldots = \vert W^{(k(n)+1)} \vert$. Тогда ясно, что $w(n) = \frac{\vert A_{1} \vert \cdot \vert A_{2} \vert}{2}  \sim \frac{n^{2}}{8}$ при $n \rightarrow \infty$. Более того, ясно, что $\vert E(W^{(i)}) \vert = \frac{1}{2} \cdot (\left[\frac{n}{4}\right]-1)\cdot w(n)$. Действительно, каждая из $w(n)$ вершин $W^{(i)}$ соединена ровно с $\left[\frac{n}{4}\right]-1$ другими вершинами из $W^{(i)}$, а сомножитель $\frac{1}{2}$ показывает, что каждое ребро было посчитано два раза. \par 
Выберем из множества $W^{(k(n)+1)}$ ровно $l(n) - k(n) \cdot \left[\frac{n}{2}\right] \cdot \left[\frac{n}{4}\right]$ вершин произвольным образом. Обозначим получившееся подмножество вершин через $U$. Ясно, что $$\vert E(U) \vert \leq \frac{1}{2} \cdot \vert U \vert \cdot \left( \left[ \frac{n}{4}\right]-1\right) \leq \frac{n^{3}}{64}.$$
Положим $$W_{n} = U \bigcup \left(\bigcup_{i=1}^{k(n)} W^{(i)} \right).$$ Тогда $$\vert W_{n}  \vert = l(n) - k(n)\cdot \left[\frac{n}{2}\right]\cdot \left[\frac{n}{4}\right] + k(n) \cdot w(n) \sim l(n) \sim k(n)\frac{n^{2}}{8}$$ при $n \rightarrow \infty$. Посчитаем мощность множества $E(W_{n})$. Обозначим $$E_{1} = \{(x,y) \in E(W_{n}) \; \vert \; \exists \; i \neq j, \; i,j \leq k(n) \; : x \in W^{(i)}, \; y \in W^{(j)}\},$$
$$E_{2} = \{ (x,y) \in E(W_{n})\; \vert \; x \in U, \; y \in W_{n} \setminus U\}.$$
 Тогда $$\vert E(W_{n}) \vert = \sum_{i=1}^{k(n)} \vert E(W^{(i)})\vert + \vert E(U) \vert+ \vert E_{1} \vert + \vert E_{2} \vert.$$ Найдем мощности множеств $E_{1}$ и $E_{2}$. Зафиксируем произвольную вершину $v \in W^{(1)}$. Обозначим $$d_{n} = \vert \{y \in W_{n} \setminus \left(W^{(1)} \cup U \right) \; \vert \; (v, y) \in E(W_{n})\} \vert.$$  Докажем, что $d_{n} = (k(n)-1) \left( \left[\frac{n}{4}\right]-2 + 2 \cdot (\vert A_{1} \vert - 1) \right)$. Действительно,  пусть $v = \{i, j, k\}, \; i \in A_{1}, \; j,k \in A_{2}$. Рассмотрим произвольную вершину $u = \{x, y, z\} \in W_{n} \setminus \left(W^{(1)} \cup U \right)$, соединенную ребром с $v$. Тогда имеют место два случая:
\begin{enumerate}
\item $x = i \text{ и } \; \vert \{j, k\} \cap \{y, z\} \vert = 0$. Cуществует $\left( k(n)-1\right)\left(\left[\frac{n}{4}\right]-2\right)$ вершин $u$, удовлетворяющих данному условию. Действительно, $ k(n)-1 $ способами можно выбрать такое натуральное $t$, что $u \in W^{(t)}$, и еще $ \left[\frac{n}{4}\right]-2 $ способами можно выбрать пару $\{y, z\}$.
\item $x \neq i \text{ и } \; \vert \{j, k\} \cap \{y, z\} \vert = 1$. Cуществует $2 \cdot \left( k(n)-1\right) \cdot \left(\vert A_{1} \vert - 1\right)$ вершин $u$, удовлетворяющих данному условию. Действительно, $k(n)-1$ способами можно выбрать такое натуральное $t$, что $u \in W^{(t)}$, еще двумя способами можно выбрать элемент, по которому пересекаются $\{j, k\}$ и $\{y, z\}$, и, наконец, $\vert A_{1} \vert - 1$ способом можно выбрать элемент $x$.
\end{enumerate}
Ясно, что $d_{n}  \sim k(n) \frac{5n}{4}$ при $n \rightarrow \infty$. Тогда в силу регулярности подграфа графа $G(n, 3, 1)$, порожденного множеством вершин $W_{n} \setminus U$, имеем
$$\vert E_{1} \vert = \frac{1}{2} \cdot d_{n} \cdot \vert W_{n} \setminus U \vert \sim \frac{5k(n)^{2}n^{3}}{64},$$
$$\vert E_{2} \vert \leq \frac{1}{2} \cdot d_{n} \cdot \vert U \vert \leq \frac{1}{2} \cdot \left[ \frac{n}{2}\right] \cdot \left[ \frac{n}{4}\right] \cdot d_{n} \leq \frac{n^{2}}{16} \cdot k(n) \cdot \left( \left[\frac{n}{4}\right]-2 + 2 \cdot (\vert A_{1} \vert - 1) \right) \sim \frac{5 k(n) n^{3}}{64}$$
при $n \rightarrow \infty$.
Итого имеем
$$\vert E(W_{n}) \vert \sim \frac{k(n)n^{3}}{64} + \frac{5k(n)^{2}n^{3}}{64} + \vert E(U) \vert + \vert E_{2} \vert \sim \frac{5k(n)^{2}n^{3}}{64} \sim \frac{5 l(n)^{2}}{\alpha_{n}}$$
при $n \rightarrow \infty$. 
Таким образом, утверждение пункта 3 доказано.

\subsection{Доказательство пункта 4}
Зафиксируем произвольную функцию $l$, удовлетворяющую условию пункта 4 теоремы, и число $n$. Положим $c_{n} = 1 - \frac{l(n)}{C_{n}^{3}}$.
Рассмотрим произвольное подмножество вершин $W \subseteq V_{n}$ мощности $l(n)$, положим $W_{1} = V_{n} \setminus W$. Ясно, что $\vert W_{1} \vert = c_{n} C_{n}^{3}$. Обозначим через $E(W_{1})$ множество ребер, концами которых являются вершины из $W_{1}$. Формально, $$E(W_{1}) = \{(x, y) \in E_{n} \; \vert \; x, y \in W_{1}\}.$$ Обозначим через $E_{1}$ множество ребер, один конец которых принадлежит множеству $W$, а другой --- множеству $W_{1}$: $$E_{1} = \{(x,y) \in E_{n} \; \vert \; x \in W, \; y \in W_{1} \}.$$ С учетом введенных обозначений мы имеем $$E(W) = E_{n} \setminus ( E(W_{1}) \sqcup E_{1} ).$$ Тогда ясно, что
$$\vert E(W) \vert = \vert E_{n} \vert - \vert E(W_{1}) \vert - \vert E_{1} \vert.$$
Оценим сверху величину $\vert E(W_{1}) \vert + \vert E_{1} \vert$. В силу регулярности графа $G(n,3,1)$ имеем $$\vert E(W_{1}) \vert + \vert E_{1} \vert \leq d_{n} \cdot \vert W_{1}\vert.$$ Действительно, каждое ребро из множеств  $E(W_{1}) \cup E_{1}$ имеет одним из своих концов вершину из  $W_{1}$. Поэтому, этих ребер не больше, чем общее число ребер, содержащих вершины из $W_{1}$. Данную оценку можно слегка уточнить. Заметим, что при таком подсчете дважды были посчитаны ребра из $E(W_{1})$.  Мощность данного множества можно оценить снизу с помощью теоремы Турана: $$\vert E(W_{1}) \vert \geq \frac{\vert W_{1} \vert^{2}}{2 \alpha_{n}} (1+o(1)).$$
Тогда $$\vert E(W_{1}) \vert + \vert E_{1} \vert \leq d_{n} \cdot \vert W_{1}\vert - \frac{\vert W_{1} \vert^{2}}{2 \alpha_{n}}(1+o(1)).$$
В итоге, суммируя все вышесказанное, имеем:
$$r(l(n)) \geq \vert E(W) \vert \geq \frac{3}{2} C_{n-3}^{2} C_{n}^{3} - d_{n} \cdot \vert W_{1} \vert + \frac{\vert W_{1} \vert^{2}}{2 \alpha_{n}}(1+o(1)) =$$ $$= \frac{3}{2} C_{n-3}^{2} C_{n}^{3} - 3 \cdot C_{n-3}^{2} \cdot \vert W_{1} \vert + \frac{\vert W_{1} \vert^{2}}{2 n} (1 + o(1)) \sim$$ $$\sim \frac{3}{2} C_{n-3}^{2} C_{n}^{3} - 3 \cdot C_{n-3}^{2} \cdot \left( \frac{c_{n} n^{3}}{6} \right) + \frac{1}{2n} \left( \frac{c_{n} n^{3}}{6}\right)^{2}(1+o(1)) \sim$$ $$\sim n^{5}\left(\frac{1}{8} - \frac{c_{n}}{4} + \frac{c_{n}^{2}}{72}\right) \; \text{при } \; n \rightarrow \infty.$$



\begin{thebibliography} {100}


\bibitem{Rai3} A.M. Raigorodskii, {\it Cliques and cycles in distance graphs and graphs of diameters}, ``Discrete Geometry and Algebraic Combinatorics'',
AMS, Contemporary Mathematics, 625 (2014), 93 - 109.

\bibitem{Rai4} A.M. Raigorodskii, {\it Coloring Distance Graphs and Graphs of Diameters}, Thirty Essays on Geometric Graph Theory, J. Pach ed., Springer, 2013, 
429 - 460.

\bibitem{Rai5} А.М. Райгородский, {\it Проблема Борсука и хроматические числа метрических пространств}, Успехи матем. наук, 56 (2001), вып. 1,
107 - 146.

\bibitem{Rai6} А.М. Райгородский, {\it О хроматических числах сфер в евклидовых пространствах}, 
Доклады РАН, 432 (2010), N2, 174 - 177.

\bibitem{Rai7} A.M. Raigorodskii, {\it On the chromatic numbers of spheres in $ {\mathbb R}^n $}, Combinbatorica, 32 (2012), N1, 111 - 123.

\bibitem{Rai8} J. Balogh, A.V. Kostochka, A.M. Raigorodskii, {\it Coloring some finite sets in $ {\mathbb R}^n $}, Discussiones Mathematicae
Graph Theory, 33 (2013), N1, 25 - 31.

\bibitem{Rai9} Л.И. Боголюбский, А.С. Гусев, М.М. Пядёркин, А.М. Райгородский, {\it Числа независимости и хроматические числа случайных подграфов в некоторых 
последовательностях графов}, Доклады РАН, 457 (2014), N4, 383 - 387. 

\bibitem{Rai10} Л.И. Боголюбский, А.С. Гусев, М.М. Пядёркин, А.М. Райгородский, {\it Числа независимости и хроматические числа случайных подграфов в некоторых 
последовательностях графов}, Матем. сборник, 2015. 

\bibitem{PA} P.K. Agarwal, J. Pach, {\it Combinatorial geometry}, John Wiley and Sons Inc., New York, 1995.

\bibitem{Szek} L.A. Sz\'ekely, {\it Erd\H{o}s on unit distances
and the Szemer\'edi--Trotter theorems}, Paul Erd\H{o}s and his Mathematics,
Bolyai Series Budapest, J. Bolyai Math. Soc., Springer, 11 (2002), 649 - 666.

\bibitem{Soi} A. Soifer, {\it The Mathematical Coloring Book}, Springer, 2009.

\bibitem{KW} V. Klee, S. Wagon, {\it Old and new unsolved problems in plane geometry and number theory}, Math. Association of America,
1991.

\bibitem{Bolt} V.G. Boltyanski, H. Martini, P.S. Soltan, 
{\it Excursions into combinatorial geometry}, Universitext, Springer, Berlin,
1997.

\bibitem{Rai1} A.M. Raigorodskii, {\it Three lectures on the Borsuk
partition problem}, London Mathematical Society Lecture Note Series,
347 (2007), 202 - 248.

\bibitem{Rai2} А.М. Райгородский, {\it Вокруг гипотезы Борсука}, Итоги науки и техники. Серия ``Современная математика'', 23 (2007), 147 - 164.

\bibitem{Ram} R.L. Graham, B.L. Rothschild, J.H. Spencer, {\it Ramsey theory}, John Wily and Sons, NY, Second Edition, 1990.

\bibitem{Nagy} Z. Nagy, {\it A certain constructive estimate of the Ramsey number}, Matematikai Lapok, 23 (1972), N 301-302, 26.

\bibitem{MS} Ф.Дж. Мак-Вильямс, Н.Дж.А. Слоэн, {\it Теория кодов, исправляющих ошибки}, М.: Радио и связь, 1979.

\bibitem{Bass} L. Bassalygo, G. Cohen, G. Z\'emor, {\it Codes with forbidden distances}, Discrete Mathematics, 213 (2000), 3 - 11.

\bibitem{Rai11} Е.Е. Демёхин, А.М. Райгородский, О.И. Рубанов, {\it Дистанционные графы, имеющие большое хроматическое число и не содержащие клик или циклов заданного 
размера}, Матем. сборник, 204 (2013), N4, 49 - 78. 

\bibitem{Rai12} А.М. Райгородский, К.А. Михайлов, {\it О числах Рамсея для полных дистанционных графов с вершинами в $ \{0,1\}^n $}, 
Матем. сборник, 200 (2009), N12, 63 - 80


\end{thebibliography}
\end{document}